\documentclass[11pt]{article}

\usepackage{amsmath,amssymb,amsthm,epsfig}

\newcommand{\vp}{\varepsilon}

\newcommand{\ovl}{\overline}

\newtheorem{thm}{Theorem}

\theoremstyle{remark}
\newtheorem{rem}{\bf{Remark}}

\begin{document}

\begin{center}
{\bf
GOOD STEIN NEIGHBORHOOD BASES AND REGULARITY\\ OF THE
$\pmb{\bar\partial}$-NEUMANN PROBLEM}
\footnote{Research supported in part by
NSF grant DMS-9801539}
\end{center}

\begin{center}
EMIL J.~STRAUBE
\end{center}
\vspace{.25in}

\centerline{1. INTRODUCTION}

The purpose of this note is to initiate a study of global regularity of the
$\bar\partial$-Neumann problem on a smooth bounded pseudoconvex domain $\Omega$
in terms of properties of Stein neighborhood bases of the closure $\ovl\Omega$ 
of the domain. If $\ovl\Omega$ is merely assumed to admit a Stein neighborhood
basis, global regularity does not follow. Indeed, the closures of the worm
domains ([\ref{D-F1}]) that wind only a little do admit Stein neighborhood
bases (see e.g. [\ref{B-F1}], Theorem~5.1), yet global regularity fails
([\ref{Ch1}, \ref{Ch2}]). On the
other hand, if $\ovl\Omega$ admits the `best possible' Stein neighborhood
basis, then the $\bar\partial$-Neumann problem is globally regular on $\Omega$.
Namely, assume that the sublevel sets of the Euclidean distance to the
boundary (outside $\ovl\Omega$), $\Omega_\vp := \{z\in
\mathbb{C}^{n}/\text{dist}(z,\ovl\Omega)< \vp\}$, are pseudoconvex for small enough
$\vp>0$. Fix such an $\vp$. Then $b\Omega$ is a level set of the distance to
the boundary of $\Omega_\vp$, hence a level set of minus the logarithm of this
distance. Because $\Omega_\vp$ is pseudoconvex, the latter function is
plurisubharmonic in $\Omega_\vp$. $\Omega$ thus admits a plurisubharmonic
defining function, and consequently the $\bar\partial$-Neumann problem is
globally regular on $\Omega$ ([\ref{B-S1}]). We show  in section~2
(Theorem~\ref{thm1}) that, more generally, if $\ovl\Omega$ admits a
`sufficiently nice' Stein neighborhood basis, then global regularity holds.
Part of the argument again exploits the plurisubharmonicity of minus the
logarithm of the boundary distance on a pseudoconvex domain. In section~3 we
briefly discuss (as it turns out to be) a generalization: \ global regularity
holds as soon as the weakly pseudoconvex directions at boundary points are
limits, from inside, of weakly pseudoconvex directions of level sets of the
boundary distance (Theorem~\ref{thm2}).

On a bounded pseudoconvex domain, the $\bar\partial$-Neumann operator $N_q$, 
$(1\le q\le n)$ is the inverse of the complex Laplacian
$\bar\partial\bar\partial^* + \bar\partial^*\bar\partial$ acting on
$(0,q)$-forms. For background on the $\bar\partial$-Neumann problem, we refer
the reader to [\ref{F-K}] and to the recent survey [\ref{B-S2}].\goodbreak
\bigskip
\goodbreak

\centerline{2. DOMAINS WITH GOOD STEIN NEIGHBORHOOD BASES}

Let $\Omega$ be a bounded smooth $(C^\infty)$ pseudoconvex domain in
$\mathbb{C}^n$. For $\vp>0$, denote by $\Omega_\vp$ the outside sublevel sets
of the Euclidean distance to $\ovl\Omega$, i.e.\ $\Omega_\vp := \{z\in
\mathbb{C}^n/\text{dist}(z,\ovl\Omega) < \vp\}$. For $\vp>0$ small enough,
$\Omega_\vp$ is also smooth. Moreover, if $n(z)$ denotes the (real) outward
unit normal to $b\Omega$, then $b\Omega_\vp = \{z+\vp n(z)/z\in b\Omega\}$.
 $W^s_{(0,q)}(\Omega)$ is the space of $(0,q)$-forms with
coefficients in $W^s(\Omega)$, with the usual inner product, where
$W^s(\Omega)$ is the $\mathcal{L}^2$-Sobolev space of order $s$. We use the
`little-o' notation in the usual sense.

\begin{thm}\label{thm1}
Let $\Omega$ be a smooth bounded pseudoconvex domain in $\mathbb{C}^n$. Assume
there  is a function $r(\vp)$ with $1-r(\vp) = o(\vp^2)$ as $\vp\to 0^+$ such
that for $\vp$ small enough, there exists a pseudoconvex domain $\widetilde
\Omega_\vp$ with $\Omega_{r(\vp)\vp} \subseteq\widetilde\Omega_\vp \subseteq
\Omega_\vp$. Then the $\bar\partial$-Neumann operators $N_q$ are continuous on
$W^s_{(0,q)}(\Omega)$, for $s\ge 0$, $1\le q\le n$.
\end{thm}

\begin{rem}\label{rem1}
The case mentioned in the introduction above, when the $\Omega_\vp$
themselves are pseudoconvex for small $\vp$, corresponds to the case
$1-r(\vp)\equiv 0$ in Theorem~\ref{thm1}. This situation was studied in
[\ref{W}]. The necessary and sufficient condition is as follows. Denote by
$H_\rho(z)$ the real Hessian of a defining function $\rho$, and by $J$ (the
matrix of) the orthogonal transformation on $\mathbb{R}^{2n}
\approx\mathbb{C}^n$ corresponding to multiplication by $\sqrt{-1}$ (on
$\mathbb{C}^n$). Assume now $\rho$ is normalized so that $|\nabla\rho| = 1$ on
$b\Omega$. Then the domains $\Omega_\vp$ are pseudoconvex for $\vp>0$ small
enough if and only if the quadratic form given by the matrix $H_\rho(z)({\bf
1} + \vp H_\rho(z))^{-1} + J^TH_\rho(z) ({\bf 1} + \vp H_\rho(z))^{-1}J$ is
positive semi-definite on the complex tangent space to $b\Omega$ at $z$ for
$z\in b\Omega$ and $\vp$ small enough. Here $J^T$ denotes the transpose of
$J$. This is proved in Proposition~1 in [\ref{W}]. Note that the leading term
(the 0-th order term) in the expansion with respect to $\vp$ is $H_\rho(z) +
J^TH_\rho(z)J$, which agrees with the Levi-form at $z$ (see e.g.\ [\ref{B}],
pp.~166--168 for a detailed comparison of the Levi-form of a hypersurface 
with its second fundamental form). \end{rem}

\begin{rem}\label{rem2}
The closure of $\Omega$ is called uniformly $H$-convex if there exists a
constant $c$, $0 < c < 1$, such that for all small enough $\vp$, there is a
pseudoconvex domain $\widetilde\Omega_\vp$ with $\Omega_{c\vp} \subseteq
\widetilde\Omega_\vp \subseteq\Omega_\vp$. This condition is useful in proving
results concerning approximation by functions holomorphic in a neighborhood of
$\ovl\Omega$, see e.g. [\ref{Ch}]. The condition used in Theorem~\ref{thm1} is
a restricted version of uniform $H$-convexity. Note that uniform $H$-convexity
is not strong enough to give global regularity: \ the worm domains that wind
only a little (that is, the parameter $r$ in [\ref{D-F1}] is sufficiently close to 1) are
uniformly $H$-convex. More specifically, for $\eta>1$ fixed, there exists
$r(\eta)$ such that if the parameter $r$ in [\ref{D-F1}] satisfies $1\le r<
r(\eta)$, then there is a defining function $\rho$ for $\Omega_r$ ($\Omega_r$
is as in [\ref{D-F1}]) such that $(\rho_+)^\eta$ is plurisubharmonic in a
neighborhood of $\ovl\Omega_r$. Here, $\rho_+(z) := \max\{\rho(z),0\}$.
This follows from the proof of Theorem~3.2 in [\ref{S}] (see also
Remark~\ref{rem2} following that proof) together with the fact that the
boundary of $\Omega_1$ is $B$-regular. (The weakly pseudoconvex boundary
points form a circle, so that the boundary is a countable union of $B$-regular
sets, hence is itself $B$-regular. See [\ref{S-2}] for details on
$B$-regularity.)
\end{rem}

The proof of Theorem~\ref{thm1} results from the existence of vector fields
that almost commute with $\bar\partial$: \ for every positive $\vp$ there
exists a vector field $X_\vp$ of type  (1,0) with coefficients in
$C^\infty(\ovl\Omega)$ such that $X\rho = 1$ on $b\Omega$ and such that when
$1\le j \le n$, the Hermitian inner product of the commutator $\left[X_\vp, 
\frac\partial{\partial\bar z_j}\right]$ with the complex normal has modulus
less than $\vp$ on the boundary. Here, $\rho$ is a defining function for
$\Omega$. It is shown in [\ref{B-S1}, \ref{B-S3}] that the existence of such a
family of vector fields implies global regularity (i.e.\ the conclusion of
Theorem~\ref{thm1}). Moreover, the proof of the lemma in [\ref{B-S1}] shows
that the crucial point is to have vectors fields with the above commutator
property, but with the commutators taken with tangential fields $\overline Z$
of type (0,1) and evaluated only at points $z\in b\Omega$ where $\overline
Z(z)$ is in the null space of the Levi form at $z$. Once  a family of vector
fields with this restricted commutator property is in hand, it can be modified
by suitable complex tangential fields to obtain a family that satisfies the
full commutator conditions (see [\ref{B-S1}], pp.~85--86 for details).

Let $\rho(z)$ be a defining function for $\Omega$ that near the boundary
equals plus or minus the boundary distance, depending on whether $z$ is
outside $\Omega$ or inside. Note that if $z\in b\Omega$, and $w\in
\mathbb{C}^n$ is contained in the complex tangent space to $b\Omega$ at $z$,
then for $|\vp|$ small, $w$ is also contained in the complex tangent space to
$b\Omega_\vp$ at $z+\vp n(z)$, where $n(z)$ denotes the outward unit normal
to $b\Omega$, and $\Omega_\vp := \{z\in \mathbb{C}^n/\rho<\vp\}$. (The real
tangent spaces at $z$ and $z+\vp n(z)$ agree, hence so do the complex tangent
spaces.) Denote by $\mathcal{L}_\rho(z)(w,\ovl w)$ the Levi form of $\rho$
(i.e.\ the complex Hessian) at $z$ applied to $w\in \mathbb{C}^n$. Let $n(z)
:= 2\left(\frac{\partial\rho}{\partial \bar z_1},\ldots,
\frac{\partial\rho}{\partial\bar z_n}\right)$, i.,e.\ $n(z)$ is the unit
normal; we also denote the associated (1,0) vector field, $2
\sum\limits^n_{j=1} \frac{\partial\rho}{\partial \bar z_j}
\frac\partial{\partial z_j}$, by $n(z)$. If $Z(z)$ is a tangential field of
type (1,0), computing the normal (1,0) component of $[n,\ovl Z]$ gives,
up to the factor 4,
\begin{align}\label{eq1}
\sum^n_{j=1} \ovl Z \left(\frac{\partial\rho}{\partial \bar z_j}\right)
\frac{\partial\rho}{\partial z_j} &= -\sum^n_{j=1}
\frac{\partial\rho}{\partial \bar z_j} \ovl Z \left(\frac{\partial
\rho}{\partial z_j}\right)\\
&= - \sum^n_{j,k=1} \frac{\partial^2\rho}{\partial z_j\partial \bar z_k}
\frac{\partial \rho}{\partial \bar z_j} \ovl \zeta_k = -\frac12
\mathcal{L}_\rho(n,\overline Z),\nonumber
\end{align}
where $Z = (\zeta_1,\ldots, \zeta_n)$. We have used that $\sum\limits^n_{j=1}
\frac{\partial \rho}{\partial z_j} \frac{\partial  \rho}{\partial \bar z_j}
\equiv 1/4$ on $b\Omega$. Consequently, by our discussion above, it suffices
to show that under the assumptions of Theorem~\ref{thm1}
\begin{equation}\label{eq2}
\mathcal{L}_\rho(z)(n,\ovl w) = 0, \qquad z\in b\Omega, w\in N_z,
\end{equation}
where $N_z$ denotes the null space of the Levi form at $z$. ($2n(z)$ will
satisfy the normalization $2n(z)(\rho)=1$.)

Consider now a weakly pseudoconvex boundary point $z_0$ and $0 \ne w\in
N_{z_0}$. We claim that 
\begin{equation}\label{eq3}
\frac{d}{d\vp} \mathcal{L}_\rho(z_0 + \vp n(z_0))(w,\ovl w)|_{\vp=0} = 0.
\end{equation}
First note that since the interior sublevel sets of $\rho$ are pseudoconvex
(by the discussion in the introduction), and $\mathcal{L}_\rho(z_0)(w,\ovl
w)=0$, this derivative is less than or equal to zero. Call its value $a$; so
$a\le 0$. The assumption that $a<0$ leads to a contradiction as follows. There
is $\vp_0>0$ such that for $0 \le \vp < \vp_0$, $\mathcal{L}_\rho(z_0+\vp
n(z_0))(w,\ovl w) \le \frac{a}2 \vp$. Thus the Levi form of $\rho$ at $z_0+\vp
n(z_0)$ has a negative eigenvalue, and this implies that analytic functions on
$\Omega_\vp$ extend past the boundary point $z_0+\vp n(z_0)$, see e.g.\
[\ref{H}], Theorem~2.6.13, where a version for extension of CR-functions is
proved. Inspection of the proof shows that there is a constant $c>0$,
independent of $\vp$ (for $0<\vp \le \vp_0$), such that one has extension in
the direction of the normal of at least $c(|a|\vp)^3$. Applying this to the
situation in Theorem~\ref{thm1}, at the point $z_0 + r(\vp)\vp n(z_0)$ (i.e.\
replacing $\vp$ by $r(\vp)\vp$), and taking into account the pseudoconvexity
of the domain $\widetilde\Omega_\vp$,  gives 
\begin{equation}\label{eq4}
r(\vp)\vp + c(|a|r(\vp)\vp)^3 \le \vp
\end{equation}
(since the extension cannot go past $b\widetilde \Omega_\vp$, hence not past
$b\Omega\vp$). (\ref{eq4}) implies
\begin{equation}\label{eq5}
c|a|^3(r(\vp))^3 \vp^2 \le 1-r(\vp),
\end{equation}
which contradicts the assumption in Theorem~\ref{thm1} that $1-r(\vp) =
o(\vp^2)$ (if $a<0$). This establishes (\ref{eq3}).

We next exploit the plurisubharmonicity of $-\log(-\rho)$ on $\Omega$, near the
boundary (where $\rho$ agrees with minus the boundary distance). The
nonnegativity of the complex Hessian of $-\log(-\rho)$ gives, after
multiplication by $-\rho$:
\begin{equation}\label{eq6}
\sum_{j,k} \left(\frac{\partial^2\rho}{\partial z_j\partial\bar z_k}(z) -
\frac1{\rho(z)} \frac{\partial\rho}{\partial z_j}(z)
\frac{\partial\rho}{\partial \bar z_k}(z)\right) w_j\ovl w_k\ge 0,\qquad w\in
\mathbb{C}^n
\end{equation}
for $z\in\Omega$ close enough to $b\Omega$. This positivity allows us to use
Cauchy-Schwarz to estimate a mixed term (as in (\ref{eq2})!) by the quadratic
terms. Moreover, the mixed term in the form in (\ref{eq6}) is the same as the
one in $\mathcal{L}_\rho$. To make this precise, denote the Hermitian form on
the left-hand side of (\ref{eq6}) by $Q_\rho(z)(w,\ovl w)$. Fix a weakly
pseudoconvex boundary point $z_0$, and let $w= (w_1,\ldots, w_n)\in N_{z_0}$.
Then for $\vp>0$ small, we have
\begin{align}
|\mathcal{L}_\rho(z_0-\vp n(z_0)) (n(z_0),w)| &= 2 \left|\sum^n_{j,k=1}
\frac{\partial^2\rho}{\partial z_j\partial \bar z_k} (z_0-\vp n(z_0))
\left(\frac{\partial\rho}{\partial \bar z_j}(z_0)\right) \ovl
w_k\right|\nonumber\\ 
\label{eq7}
&= 2|Q_\rho(z_0-\vp n(z_0)) (n(z_0), \overline w)|\\
&\le 2 Q_\rho(z_0-\vp n(z_0)) (n(z_0), \ovl{n(z_0)})^{1/2} Q_\rho(z_0-\vp
n(z_0)) (w,\ovl w)^{1/2}\nonumber\\
&= 2 Q_\rho(z_0-\vp n(z_0)) (n(z_0), \ovl{n(z_0)}) ^{1/2}
\mathcal{L}_\rho(z_0-\vp n(z_0)) (w,\ovl w)^{1/2}.\nonumber
\end{align}
Note that the normal to $b\Omega_{-\vp}$ at $z_0-\vp n(z_0)$ equals $n(z_0)$,
and that $w$ also belongs to the complex tangent space to $b\Omega_{-\vp}$ at
$z_0-\vp n(z_0)$, so that $\sum\limits^n_{k=1} \frac{\partial\rho}{\partial
\bar z_k} (z_0-\vp n(z_0)) \ovl w_k =0$. Consequently, $\mathcal{L}_\rho(z_0-
\vp n(z_0))(\cdot, \ovl w) = Q_\rho (z_0-\vp n(z_0))(\cdot, \ovl w)$. The
quadratic terms in the last term in (\ref{eq7}) are now easily estimated.
First, $Q_\rho(z_0-\vp n(z_0)) (n(z_0), \ovl{n(z_0)})$ is dominated by
$1/\vp$, just by the definition of $Q_\rho$. Second, $\mathcal{L}_\rho(z_0-\vp
n(z_0))(w,\ovl w)$ is dominated by $\vp^2$, in view of (\ref{eq3}) and the
fact that $w\in N_{z_0}$ (i.e.\ $\mathcal{L}_\rho(z_0)(w,\ovl w) = 0)$.
Inserting these  two estimates into (\ref{eq7}) and letting $\vp\to 0^+$
yields (\ref{eq2}). This completes the proof of Theorem~\ref{thm1}.

\smallskip
\begin{rem}\label{rem3}
The argument in the last part of the proof of Theorem~\ref{thm1} that exploits
the plurisubharmonicity of $-\log(-\rho)$ is of interest in part because of
its potentially broader scope, but in the situation of Theorem~\ref{thm1}, one
can exploit the work in [\ref{W}] to obtain an alternative proof (as well as
further information). It is shown in [\ref{W}] that $\mathcal{L}_\rho(z_0+\vp
n(z_0))(w,\ovl w)$ can be conveniently expressed in terms of $H_\rho(z_0)$,
the {\em real\/} Hessian of $\rho$ at $z_0$. In fact, the expression is the
one quoted in Remark~1 above. (One can always express the Levi form in terms of
the second fundamental form, see [\ref{B}], pp.~166--168. The point of
[\ref{W}] is that in the case of the level surfaces of the boundary distance,
the (matrix corresponding to the) second fundamental form of the level surface
is easily expressed in terms of the (matrix corresponding to the) second
fundamental form of the boundary itself). From this expression, the derivative
 in (\ref{eq3}) is computed to be  
\begin{equation}\label{eq8}
\frac{d}{d\vp}\mathcal{L}_\rho (z_0+\vp n(z_0))
(w,\ovl w)|_{\vp=0} = -(\|H_\rho(z_0) \omega\|^2 + \|H_\rho(z_0) J\omega\|^2),
\end{equation}
where $\omega$ denotes the vector in $\mathbb{R}^{2n}$ corresponding to $w\in
\mathbb{C}^n$ under the usual identification. Thus if this derivative
vanishes, we obtain
\begin{equation}\label{eq9}
H_\rho(z_0) \omega = H_\rho(z_0)J\omega = 0,\qquad \omega\in N_{z_0}.
\end{equation}
That is, (real) derivatives of the normal in directions in the null space of
the Levi form vanish. By the discussion above (see in particular (\ref{eq1})),
this suffices to prove the theorem. This approach gives some additional
information. In particular, in the situation of Theorem~\ref{thm1}, if $z_0\in
b\Omega$ and $w\in N_{z_0}$, then $\mathcal{L}_\rho(z_0+\vp n(z_0))(w,\ovl w)
\equiv 0$ for $\vp$ close enough to 0. It is instructive to compare this with
the necessary and sufficient condition for pseudoconvexity of $\Omega_\vp$
discussed in Remark~\ref{rem1} above; compare also the discussion on p.~401 in
[\ref{W}], in particular Corollary~2.
\end{rem}
\medskip

\goodbreak
\centerline{3. A GENERALIZATION}

The arguments in the proof of Theorem~\ref{thm1} show that global regularity
holds if (\ref{eq3}) holds, but (see Remark~3 above), if (\ref{eq3}) holds,
then actually $\mathcal{L}_\rho(z_0+\vp n(z_0))(w,\bar w) \equiv 0$ for $|\vp|$
small and $w\in N_{z_0}$. In particular, the weakly pseudoconvex directions at
boundary points are limits, along the normal, of weakly pseudoconvex
directions at points of interior level sets of the boundary distance. The
following is thus a generalization (also of Theorem~\ref{thm1}).
We keep the notation from above, so that $\rho$ is still the defining
function given by $+$ or $-$ the boundary distance (near $b\Omega$).

\begin{thm}\label{thm2}
Let $\Omega$ be a smooth bounded pseudoconvex domain in $\mathbb{C}^n$ with the
property that for each pair $(z,w)\in b\Omega\times \mathbb{C}^n$ with $w\in
N_z$, there is a sequence $\{(z_n,w_n)\}^\infty_{n=1}$ in $\Omega\times
\mathbb{C}^n$ such that $(z_n,w_n)\to (z,w)$ as $n\to\infty$, and for all
$n,w_n$ is in the null space of the Levi-form at $z_n$ of the level set of
$\rho$  through $z_n$. Then the $\bar\partial$-Neumann operators $N_q$, $1\le
q\le n$, are continuous on $W^s_{(0,q)}(\Omega)$, for $s\ge 0$.
\end{thm}

The proof of Theorem~\ref{thm2} is the same as that of Theorem~\ref{thm1},
with the obvious modification in (\ref{eq7})
\begin{equation}\label{eq10}
\mathcal{L}_\rho(z_n)(n(z_n), \ovl w_n) = 2Q_\rho(z_n)(n(z_n),\ovl w_n) = 0.
\end{equation}
In the last equality in (\ref{eq10}), we have used that $w_n$ is also a null
direction of $Q_\rho$ (since $Q_\rho(z_n)(w_n,\ovl w_n) =
\mathcal{L}_\rho(z_n)(w_n,\ovl w_n))$ and again that $Q_\rho$ is positive
semi-definite. Letting $n$ tend to $\infty$ in (\ref{eq10}) gives (\ref{eq2}).

\smallskip
\begin{rem}\label{rem4}
 In view of the consequences,
it would be of interest to have a better understanding of the condition in 
Theorem~\ref{thm2}. An
easy observation is that it is satisfied by convex domains (for which global
regularity is of course known, [\ref{B-S1}]). In fact, convex domains satisfy
the stronger condition (discussed in the introduction) that the outside
sublevel sets of the Euclidean distance are pseudoconvex, see [\ref{W}],
Corollary~2 and the discussion following it.
\end{rem}
\bigskip

\baselineskip = 12pt

\newpage
\centerline{REFERENCES}

\begin{enumerate}
\item\label{B-F1}
 Eric Bedford and John Erik Fornaess: \ Domains with pseudoconvex neighborhood
systems, Inventiones Math. {\bf 47} (1978), 1--27.

\item\label{B-S1}
Harold P.\ Boas and Emil J.\ Straube: \ Sobolev estimates for the
$\bar\partial$-Neumann operator on domains in $\mathbb{C}^n$ admitting a
defining function that is plurisubharmonic on the boundary, Math. Z. {\bf 206}
(1991), 81--88.

\item\label{B-S3}
Harold P.\ Boas and Emil J.\ Straube: \ De Rham cohomology of manifolds
containing the points of infinite type, and Sobolev estimates for the
$\bar\partial$-Neumann problem, Journal Geom.\ Anal. {\bf 3}, Nr.~3 (1993),
225--235.

\item\label{B-S2}
Harold P.\ Boas and Emil J.\ Straube: \ Global regularity of the
$\bar\partial$-Neumann problem: \ a survey of the $\mathcal{L}^2$-Sobolev
theory, Several Complex Variables, M.~Schneider and Y-T.\ Siu, editors, MSRI
Publications {\bf 37}, 79--111, Cambridge Univ.\ Press, 1999.

\item\label{B}
Albert Boggess: \ CR Manifolds and the Tangential Cauchy-Riemann Complex, CRC
Press, 1991.

\item\label{Ch1}
Michael Christ: \ Global $C^\infty$ irregularity of the $\bar\partial$-Neumann
problem for worm domains, J.\ Amer.\ Math.\ Soc. {\bf 9}, Nr. 4 (1996),
1171--1185.

\item\label{Ch2}
Michael Christ: \ Remarks on global irregularity in the $\bar\partial$-Neumann
problem, Several Complex Variables,  M.~Schneider and Y.-T.\ Siu, editors,
MSRI Publications {\bf 37}, 161--198, Cambridge Univ.\ Press, 1999.

\item\label{Ch}
E.M.\ \u Cirka: \ Approximation by holomorphic functions on smooth manifolds in
$\mathbb{C}^n$, Mat.\ Sb. {\bf 78} (1969), 101--123, Math.\ USSR-Sb. {\bf 7}
(1969), 95--113.

\item\label{D-F1}
Klas Diederich and John Erik Fornaess: \ Pseudoconvex domains: \ an example
with nontrivial Nebenh\"ulle, Math.\ Ann. {\bf 225} (1977), 275--292.

\item\label{F-K}
G.B.\ Folland and J.J.\ Kohn: \ The Neumann Problem for the Cauchy-Riemann
Complex, Annals of Math.\ Studies {\bf 75}, Princeton Univ.\ Press, 1972.

\item\label{H}
Lars H\"ormander: \ Complex Analysis in Several Variables, $2^{\rm nd}$
edition, North Holland Publishing Company, 1973.

\item\label{S-2}
Nessim Sibony: \ Une classe de domaines pseudoconvexes, Duke Math.\ J. {\bf
55}, Nr.~2 (1987), 299--319.

\item\label{S}
Nessim Sibony: \ Some aspects of weakly pseudoconvex domains, Several Complex
Variables and Complex Geometry, E.\ Bedford, J.\ D'Angelo, R.\ Greene, and S.\
Krantz, editors, Proc.\ Symp.\ Pure Math., {\bf 52}, part 1, 199--231,
Amer.\ Math.\ Soc., 1991.

\item\label{W}
Barnet M.\ Weinstock: \ Some conditions for uniform $H$-convexity, Illinois
J.\ Math. {\bf 19} (1975), 400--404.
\end{enumerate}

\vskip.5in

{\parindent=0pt\obeylines
Department of Mathematics
Texas A\&M University
College Station, TX \ 77843-3368
straube@math.tamu.edu}

\end{document}